\documentclass{article}
\usepackage{amsthm,amsfonts, amsbsy, amssymb,amsmath,graphicx}
\usepackage{epsfig}
\usepackage{graphicx}

\begin{document}

\newcommand{\fh}{\frac{1}{2}}
\newcommand{\beq}{\begin{equation}}
\newcommand{\eeq}{\end{equation}}

\date{}

\title{\bf Particle Topology, Braids and Braided Belts}

\author{Sundance Bilson-Thompson\\Perimeter Institute for Theoretical Physics\\31 Caroline Street North\\
Waterloo, Ontario N2L 2Y5,Canada\\ \\
Jonathan Hackett, \\Perimeter Institute for Theoretical Physics\\31 Caroline Street North\\
Waterloo, Ontario N2L2Y5,Canada and\\ Department of Physics\\ Unversity of Waterloo\\ Waterloo, Ontario N2J 2W9,Canada\\ \\
Louis H. Kauffman\\ Department of Mathematics, Statistics \\ and Computer Science (m/c
249)    \\ 851 South Morgan Street   \\ University of Illinois at Chicago\\
Chicago, Illinois 60607-7045, USA\\ \\
}

 \maketitle
  
 \thispagestyle{empty}
 
 \subsection*{\centering Abstract}

{\em Recent work \cite{Sundance,Sundance1} suggests that topological features of certain quantum gravity theories can be interpreted as particles, matching the known fermions and bosons of the first generation in the Standard Model. This is achieved by identifying topological structures with elements of the framed Artin braid group on three strands, and demonstrating a correspondence between the invariants used to characterise these braids (a braid is a set of non-intersecting curves, that connect one set of $N$ points with another set of $N$ points), and quantities like electric charge, colour charge, and so on \cite{Sundance,Hackett2}. In this paper we show how to manipulate a modified form of framed braids to yield an invariant standard form for sets of isomorphic braids, characterised by a vector of real numbers. This will serve as a basis for more complete discussions of quantum numbers in future work.} 
 
\section{Introduction} It has long been hoped that the fundamental particles, of which all matter in the Universe is composed, would turn out to be topological structures of some type.
This prospect is appealing because if correct it would indicate that the permanence of matter and its properties (obeying principles like the conservation of electric charge) could be explained in terms of different topological classes, which are isotopically inequivalent. Furthermore it would explain the number and type of different particles (electrons, neutrinos, quarks) in terms of a simple counting exercise, much as the periodic table of the elements can be explained by counting the number of protons in atomic nuclei. In \cite{Sundance} it was proposed that certain hypothetical models of particle substructure could be reformulated in terms of the framed Artin braid group on three strands. This idea was soon adapted to interpret extra topological degrees of freedom occurring in certain theories of quantum gravity \cite{Sundance1}, raising the exciting possibility that matter may be an emergent feature of quantum gravity theories. In order to test this idea, it is necessary to understand how such braided structures may be classified into isotopically inequivalent classes, what quantities characterise these classes, and how these may correspond with the quantum numbers of the standard model. 
\bigbreak
In this paper we present a general method for reducing framed braids on three strands, carried on a closed surface, to a simplified form, in which all the crossings have been removed, and only twists remain on the strands. Each braid can then be classified into an equivalence class, each such class being specified by a triplet $[a,\,b,\,c]$ of multiples of half-integers. In section~\ref{sec:artin} we describe the Artin braid group, the algebra of braids, and the generators of the braid group on $n$ strands. In section~\ref{sec:belts} we extend the braid group to framed braids, and introduce the physical example of braided belts (these being framed braids carried on a closed surface, the ``belt''), which provide a useful mental image to complement the mathematical discussion. We demonstrate how braided belts on three strands may be reduced to a simplified form (without crossings), classified by a triplet of real numbers. In section~\ref{sec:algebra} we reiterate several of the results from the previous section in a mathematical formalism utilising permutation matrices. We conclude with a demonstration of the classification of the quarks and leptons of the first generation of the standard model in section~\ref{sec:concl}.

\section{The Artin Braid Group}
\label{sec:artin}

\begin{figure}[h]
     \begin{center}
     {\includegraphics[height=40mm,angle=0]{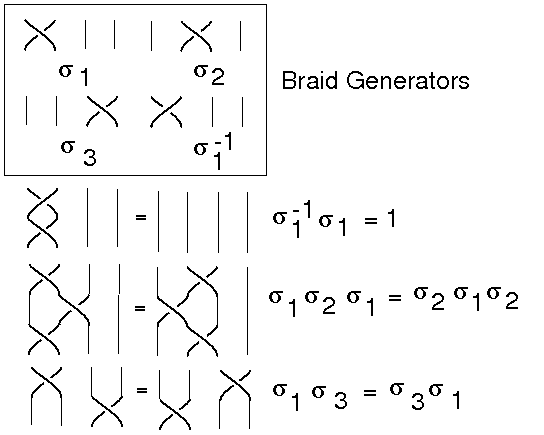}}
     \caption{Braid Generators and Relations }
     \label{fig:generators}
\end{center}
\end{figure}

Let $B_{n}$ denote the Artin braid group on $n$ strands \cite{Sossinsky}.
We recall here that $B_{n}$ is generated by elementary braids $\{\sigma_1,\ldots,\sigma_{n-1}\}$
with relations 

\begin{enumerate}
\item $\sigma_{i} \sigma_{j} = \sigma_{j} \sigma_{i}$ for $|i-j| > 1$, 
\item $\sigma_{i} \sigma_{i+1} \sigma_{i} = \sigma_{i+1} \sigma_{i} \sigma_{i+1}$ for $i= 1, \cdots , n-2.$
\end{enumerate}

\begin{figure}
     \begin{center}
     \begin{tabular}{c}
     \includegraphics[width=7cm]{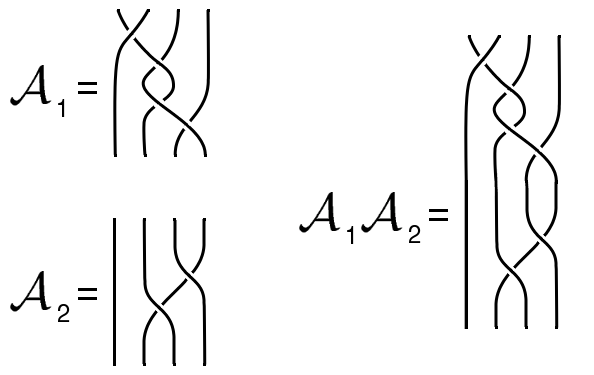}
     \end{tabular}
     \caption{The product of two braids, ${\cal A}_1$ and ${\cal A}_2$.}
     \label{fig:braidproduct}
\end{center}
\end{figure}

\noindent See Fig.~\ref{fig:generators} for an illustration of the elementary braids and their relations. Note that the braid group has a diagrammatic topological interpretation, as follows. Consider a set $P$ consisting of $n$ distinct points, $\{p_1,\,p_2,\ldots,\,p_n\}$ sharing the same $z$-coordinate (say $z_P$), and a second set $Q$ of $n$ distinct points, $\{q_1,\,q_2,\ldots,\,q_n\}$ sharing a different $z$-coordinate (say $z_Q$). Then a braid is a set of $n$ non-intersecting strands that lead from the first set of $n$ points to the second set, such that the $z$-coordinate of a point along any strand changes {\it monotonically} from $z_P$ to $z_Q$ as we move from a point in $P$, to the corresponding point in $Q$. In other words, strands cannot ``loop back'' on (and be tied in knots around) themselves. The braid generators $\sigma_i$ are represented by diagrams where the $i^\mathrm{th}$ and $(i + 1)^\mathrm{th}$ strands wind around one another by a single half-twist (the sense of this winding is shown in Fig.~\ref{fig:generators}) and all other strands drop straight to the bottom. Braids are diagrammed vertically as in Fig.~\ref{fig:generators}, and the products are taken in order from top to bottom. The product of two braid diagrams is accomplished by adjoining the tops of the strands of the second braid to the bottoms of the strands of the first braid, or in other words, identifying the $n$ points in the first set of the second braid with the points in the second set of the first braid. An example of the product of two braids is illustrated in Fig.~\ref{fig:braidproduct}
\bigbreak 

Any braid on $n$ strands can be written as a product of the generators $\sigma_1,\ldots,\sigma_{n-1}$ and their inverses. The sequence of $\sigma$ factors defining a braid is called its {\em braid word}. We shall sometimes use the terms ``braid'' and ``braid word'' interchangeably, for brevity.
\bigbreak

In Fig.~\ref{fig:generators} we illustrate the example of the four-stranded braid group $B_4$. In that figure the three braid generators of $B_4$ are shown, and then the inverse of the first generator is drawn. Following this, one sees the identities $\sigma_{1} \sigma_{1}^{-1} = I$ 
(where the identity element in $B_{4}$ consists of four vertical strands), 
$\sigma_{1} \sigma_{2} \sigma_{1} = \sigma_{2} \sigma_{1} \sigma_{2},$ and finally
$\sigma_1 \sigma_3 = \sigma_3 \sigma_1.$ 
\bigbreak

We also observe that the generators induce permutations of the strand ordering. The generator $\sigma_1$ induces the permutation $P_{1,2}$ (that is, it swaps the $1^\mathrm{st}$ and $2^\mathrm{nd}$ strands), the generator $\sigma_2$ induces the permutation $P_{2,3}$ and in general the generator $\sigma_i$ induces the permutation $P_{i,(i+1)}=P_{(i+1),i}$. Notice also that the same permutation is induced by a generator or its inverse, $\sigma^{-1}_i$. Therefore the generators contain more information than the permutations - in particular the direction of the crossing is specified by the generators (as shown in Fig.~\ref{fig:generators}), the inverse of a generator inducing a crossing which is opposite in direction to the crossing induced by the generator itself.
\bigbreak

For further discussions of braids (and knots) the reader is directed to references \cite{Sossinsky,AT,Birman,Conway,KP,WIT}.
\bigbreak

\section{Braided Belts}
\label{sec:belts}
\begin{figure}
     \begin{center}
     \begin{tabular}{c}
     \includegraphics[width=10cm]{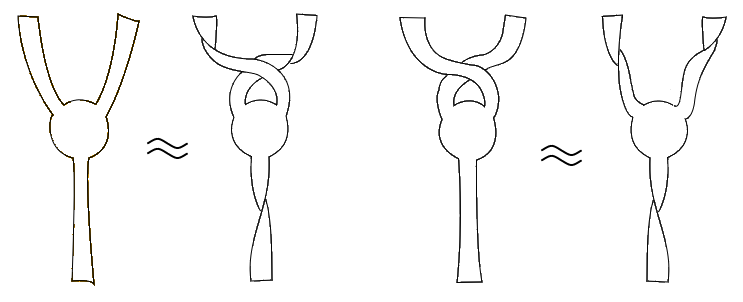}
     \end{tabular}
     \caption{Turning over a node induces crossings and twists in the ``legs" of that node (left). A trinion flip may be used to eliminate crossing while creating twists (right).}
     \label{fig:twistY}
\end{center}
\end{figure}

\begin{figure}
     \begin{center}
     \begin{tabular}{c}
     \includegraphics[width=7cm]{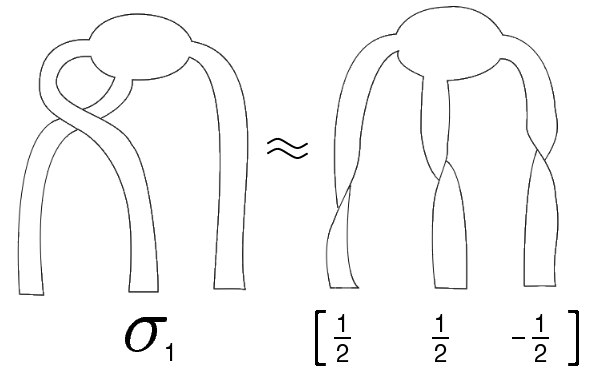}
     \end{tabular}
     \caption{Trinion bent to form a generator of the braid group}
     \label{fig:trinionGen}
\end{center}
\end{figure}

We now turn to studying surfaces that we shall call ``braided belts". As Figs~\ref{fig:twistY}-\ref{fig:belt} indicate, the braiding we are considering is a braiding 
of ribbon strands. This requires a natural generalization of the braid group to the 
{\it framed braid group} where each strand of the braid has been replaced by a ribbon, which can carry twist. We count the twist on each strand or ribbon in terms of the number of {\it half-twists} in the ribbon (i.e. twists through $180^\circ$). It is therefore convenient to define a standardised form for framed braids in which all the twist is isotoped to the top of the braid. Then we can write $[r,s,t]{\cal B}$ where ${\cal B}$ is an ordinary braid word and $[r,s,t]$ is a triple of multiples of half-integers which catalogue the twists in the ribbons. We shall call this triple of numbers the {\em twist-word}. Thus a framed braid on three strands is completely specified by the twist word and the braid word (this can of course be easily generalised to framed braids on $n$ strands).
\bigbreak 

At this point it is necessary to draw a distinction between framed braids, and what we call braided belts. Framed braids are simply those braids which consist of $n$ ribbons connecting one set of $n$ points with a common $z$-coordinate, to another set of $n$ points also with a common $z$-coordinate, as described above. In braided belts the sets of $n$ points at both ends of the ribbon strands are replaced by disks, or framed nodes. The union of the strands and disks is a closed surface. The disks are thus nodes with $n$ ``legs'' emerging from them, each leg being the terminus of a ribbon strand. We furthermore relax the condition that strands may not ``loop back'' on themselves. The ability to deform the braided belt by flipping the node at one end over (effectively feeding it through the strands), in order to undo crossings and hence simplify the braid structure of the belt, is essential to the results we will discuss below.
\bigbreak

On the left of Fig.~\ref{fig:twistY} we illustrate a 
{\bf Y}-shaped ribbon - that is, a node with three ribbons or ``legs'' attached to it. Such a 
triple of ribbons meeting at a node has previously been referred to as a ``trinion'' \cite{Sundance1}, and we shall adopt this terminology. A trinion may be converted into a structure with both crossings and twists, by keeping the ends fixed and flipping over the node in the middle, as illustrated in Fig.~\ref{fig:twistY} (we shall refer to this process of flipping over a node while keeping the ends of the legs fixed as a ``trinion flip'', or ``trip'' for short). Conversely, on the right of Fig.~\ref{fig:twistY} we show how a trinion with untwisted ribbons, but whose upper ribbons are crossed, can be converted into a trinion with uncrossed ribbons and oppositely-directed half-twists in the upper and lower ribbons by performing an appropriate trinion flip (in the illustration, a negative half-twist in the lower ribbon of the trinion and positive half-twists in the upper ribbons). In Fig.~\ref{fig:trinionGen}, we show the same process performed on a trinion whose (crossed) upper ribbons have been bent downwards to lie besides and to the left of the (initially) lower ribbon. This configuration is nothing other than a framed 3-braid corresponding to the generator $\sigma_1$ (with the extra detail that the tops of all three strands are joined at a node). Keeping the ends of the ribbons fixed as before and flipping over the node so as to remove the crossings now results in three unbraided (i.e. trivially braided) strands, with a positive half-twist on the leftmost strand, a  positive half-twist on the middle strand, and a negative half-twist on the rightmost strand. Hence the associated twist-word is $[\fh,\fh,-\fh]$ This illustrates that in the case of braids on three strands, each of the crossing generators can be isotoped to uncrossed strands bearing half-integer twists. By variously bending the top two ribbons down to the right of the bottom ribbon, and/or taking mirror images, and performing the appropriate trinion flips we can determine that the generators may be exchanged for twists according to the pattern;
\begin{equation}
\begin{tabular}{lll}
$\sigma_1$      & $\rightarrow$ & $\left[\fh, \:\fh,\:-\fh\right]$ \vspace{1mm} \\
$\sigma_1^{-1}$ & $\rightarrow$ & $\left[-\fh,\:-\fh,\:\fh\right]$ \vspace{1mm} \\
$\sigma_2$      & $\rightarrow$ & $\left[-\fh,\:\fh, \:\fh\right]$ \vspace{1mm} \\
$\sigma_2^{-1}$ & $\rightarrow$ & $\left[\fh,\:-\fh,\:-\fh\right]$
\end{tabular}
\label{eq:sigmaToTwist}
\end{equation}
All braids on three strands can be built up as products of these generators.
 
\subsection{Framed braid multiplication}
As mentioned above, a braid with several crossings is specified by its braid word, which corresponds with several generators multiplied together. Unframed braids are multiplied together by joining the tops of the strands of the second braid to the bottoms of the strands of the first braid. Framed braids are multiplied in an analogous way, however there are extra complications introduced by the presence of twists. Therefore when forming the product of two (or more) framed braids, we shall firstly join the bottom of the ribbons in the first braid to the tops of the ribbons in the second braid, and so on (in the case of braided belts, we must eliminate the intervening nodes first, while retaining the nodes at the ``top'' of the first belt in the product, and the ``bottom'' of the last belt). Secondly we shall isotop the twists from each of the component braids upwards to render the product into standard form (as described above). They will get permuted in the process, according to the permutation $P_{\cal B}$ associated with the braid ${\cal B}$ that they pass through. Given a permutation $P_{\cal B}$ and a twist-word $[x,y,z]$, we shall write the permuted twist-word as $P_{\cal B}\left([x,y,z]\right)$. Thus given two braids $[r,s,t]{\cal B}_1$ and $[u,v,w]{\cal B}_2$, we can form their product by joining strands to form 
$[r,s,t]{\cal B}_1[u,v,w]{\cal B}_2$, and then move the twists $[u,v,w]$ upward along the strands of the braid ${\cal B}_1$.  Thus 
\begin{equation}
[r,s,t]{\cal B}_1[u,v,w]{\cal B}_2 = [r,s,t]P_{{\cal B}_1}\left([u,v,w]\right){\cal B}_1{\cal B}_2
\label{eq:wordpermute}
\end{equation}
where $[r,s,t][x,y,z] = [r+x,s+y,t+z]$ and ${\cal B}_1{\cal B}_2$ denotes the usual product of 
braid words. In general the twists will be permuted by all the braids they pass through. For example, remembering that $P_{\sigma_i} = P_{i,(i+1)}$, 
\begin{eqnarray} 
[r,s,t]\sigma_1\sigma_2[x,y,z] & = & [r,s,t]\sigma_1[x,z,y]\sigma_2 \nonumber \\
                               & = & [r,s,t][z,x,y]\sigma_1\sigma_2 \nonumber \\
			       & = & [r+z,s+x,t+y]\sigma_1\sigma_2.
\label{eq:twistpermut}
\end{eqnarray}
It should also be clear to the reader that since crossings can be exchanged for twists, in the case of braided belts - which we shall be discussing exclusively from this point onwards, we may go a step further and entirely eliminate the crossings from a 3-braid. When we do so we uncross the strands (hence permuting them by the permutation associated with the crossing being eliminated) and introduce the twists indicated in eqn.~(\ref{eq:sigmaToTwist}). In general, this means that we iterate the process
\begin{eqnarray}
[a_1, a_2, a_3][b_1, b_2, b_3]\sigma_i\sigma_j\ldots\sigma_m 
\rightarrow [a_1+b_1, a_2+b_2, a_3+b_3]\sigma_i\sigma_j\ldots\sigma_m \nonumber & & \\
\:\:\:\:\:\: \rightarrow P_{\sigma_i}([a_1+b_1, a_2+b_2, a_3+b_3])[x,y,z]\sigma_j\ldots\sigma_m & & 
\label{eq:sigmaTwistIteration}
\end{eqnarray}
where $[x,y,z]$ is the twist-word associated to $\sigma_i$ (as listed in Eqn.~\ref{eq:sigmaToTwist}, when $i$ is specified). We iterate this procedure until the braid word becomes the identity. Hence continuing the example above, from Eqn.~(\ref{eq:twistpermut}), 
\begin{eqnarray}
[r+z,s+x,t+y]\sigma_1\sigma_2 & \rightarrow & P_{\sigma_1}([r+z,s+x,t+y])\left[\fh,\fh,-\fh\right]\sigma_2\nonumber \\
	& \rightarrow & [s+x,r+z,t+y]\left[\fh,\fh,-\fh\right]\sigma_2 \nonumber \\
	& \rightarrow & \left[s+x+\fh, r+z+\fh, t+y-\fh\right]\sigma_2 \nonumber \\
	& \rightarrow & P_{\sigma_2}\left(\left[s+x+\fh, r+z+\fh, t+y-\fh\right]\right) \nonumber \\
	&   & \:\:\:\:\:\:\:\:\:\:\:\:\:\:\:\:\:\:\:\:\:\:\:\:\times \left[-\fh,\fh,\fh\right]\nonumber \\
	& \rightarrow & \left[s+x+\fh,t+y-\fh,r+z+\fh\right]\nonumber \\
	&   & \:\:\:\:\:\:\:\:\:\:\:\:\:\:\:\:\:\:\:\:\:\:\:\:\times \left[-\fh,\fh,\fh\right]\nonumber \\
	& \rightarrow & [s+x,t+y,r+z+1].
\end{eqnarray}
We shall refer to the form of a braid in which all the crossings have been exchanged for twists as the {\em pure twist form}. The list of three numbers which characterise the twists on the strands in the pure twist form will be referred to as the {\em pure twist-word}. The pure twist-word is of interest because it is a topological invariant (since it is obtained when a  braid is reduced to a particularly simple form i.e. all crossings removed).
\bigbreak

\subsection{Making 3-belts}
\begin{figure}
     \begin{center}
     \begin{tabular}{c}
     \includegraphics[width=5cm]{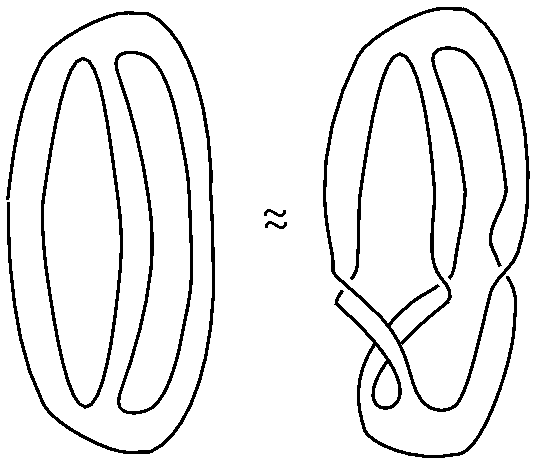}
     \end{tabular}
     \caption{A Trip performed on a closed 3-Belt}
     \label{fig:closed}
\end{center}
\end{figure}

Consider a braided belt (or framed braid) on three strands. In the particular (trivial) case where the strands do not cross each other, the associated braid word is clearly the identity, ${I}$. In the case where the strands are untwisted, the associated twist-word is also the identity. Such an untwisted trivial braid - the identity braid on three strands - can be made by cutting two parallel slits in a strip of leather, as shown on the left of Fig.~\ref{fig:closed}. The resulting surface is topologically equivalent to three parallel strips capped at the top and bottom by an attached disk. We will refer to three strands (not necessarily unbraided) attached to disks in this manner as a {\em 3-belt}.
In Fig.~\ref{fig:closed}, we show the consequence of trinion flips in the making of a braided leather belt. 
\bigbreak

\begin{figure}
     \begin{center}
     \begin{tabular}{c}
     \includegraphics[width=7cm]{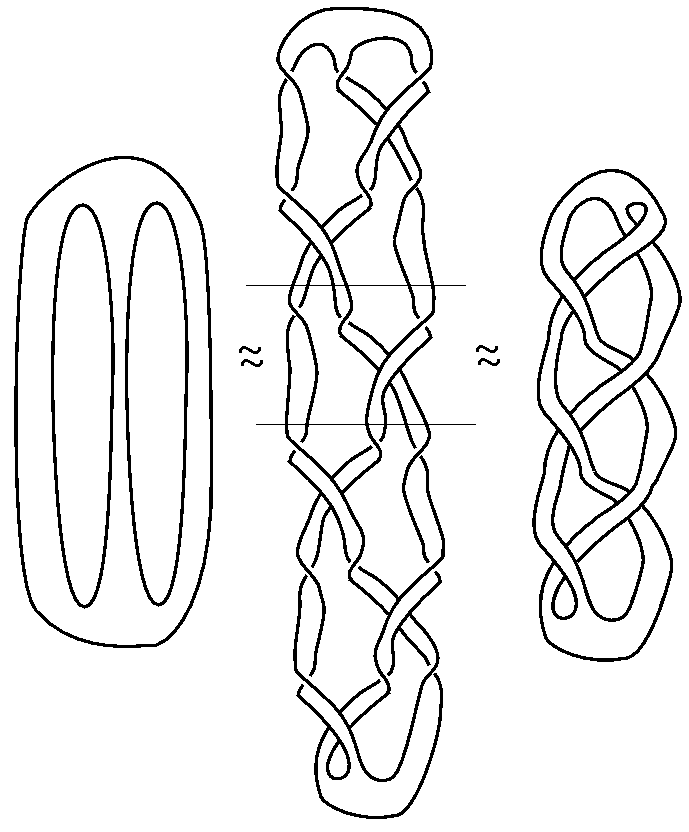}
     \end{tabular}
     \caption{Braiding a Belt}
     \label{fig:belt}
\end{center}
\end{figure}

Recall that $\sigma_1 \equiv [\fh,\,\fh,\,-\fh]$, and so we may write 
${I} = [-\fh,\,-\fh,\,\fh]\sigma_1$. The right-hand side of Fig.~\ref{fig:closed} illustrates that the strands are now crossed and twisted, but still the 3-belt we have obtained is isotopic to the trivial 3-belt. Figure~\ref{fig:belt} shows the result of six consecutive repetitions of this process (alternately to the first two strands and the last two strands) of a 3-belt. The reader will note by direct observation that along each of the three ribbon strands, all the twists cancel. Thus when we isotop all the twists to the top of the braid we obtain a flat (untwisted) braided belt, with braid-word $(\sigma^{-1}_2\sigma_1)^3$. Iteration of the procedure that yields Fig.~\ref{fig:belt} is actually used by belt-makers. 
\bigbreak

Since the same physical structure can be isotoped to have the braid word $I$ or $(\sigma^{-1}_2\sigma_1)^3$, it is clear that the braid word is not a topological invariant. However, as noted above, the pure twist-word {\it is} a topological invariant. If any two braids $[a,b,c]{\cal B}_1$ and $[x,y,z]{\cal B}_2$ have the same pure twist-word, then they are isotopic. For the remainder of this paper we shall be mostly interested in classifying braided, twisted 3-belts by their pure twist numbers, rather than inducing twists and crossings on an initially trivial 3-belt. We are therefore primarily interested in the procedure illustrated in Eqn.~(\ref{eq:sigmaTwistIteration}). We now apply this procedure to the braid word $(\sigma^{-1}_2\sigma_1)^3$, and confirm that its pure twist-word is $I=[0,0,0]$. This corresponds with reversing the procedure (used by belt-makers) described in the previous paragraph. 
\bigbreak

Firstly we consider the braid word $\sigma^{-1}_2\sigma_1$
\begin{eqnarray}
\sigma_{2}^{-1}\sigma_{1} & \rightarrow & \left[\fh,-\fh,-\fh\right]\sigma_{1} \nonumber \\
 & \rightarrow & P_{1,2}\left(\left[\fh,-\fh,-\fh\right]\right)\left[\fh,\fh,-\fh\right] \nonumber \\
 & \rightarrow & \left[-\fh,\fh,-\fh\right]\left[\fh,\fh,-\fh\right] \nonumber \\
 & \rightarrow & [0,1,-1]
\end{eqnarray}
To consider the full braid with six crossings we need to multiply this result with itself three times, but also keep track of the permutations induced by the braid word $\sigma^{-1}_2\sigma_1$. This permutation will be $P_{1,2}(P_{2,3})=P_{(123)}$, that is, the cyclic permutation which sends $[a,b,c]\rightarrow[c,a,b]$. Hence,
\begin{eqnarray}
(\sigma_{2}^{-1}\sigma_{1})^3 & \rightarrow & [0,1,-1] (\sigma_{2}^{-1}\sigma_{1})^2 \nonumber \\
 & \rightarrow & P_{(123)}([0,1,-1])[0,1,-1]\sigma_{2}^{-1}\sigma_{1} \nonumber \\
 & \rightarrow & [-1,0,1][0,1,-1]\sigma_{2}^{-1}\sigma_{1} \nonumber \\
 & \rightarrow & [-1,1,0]\sigma_{2}^{-1}\sigma_{1} \nonumber \\
 & \rightarrow & P_{(123)}([-1,1,0])[0,1,-1] \nonumber \\
 & \rightarrow & [0,-1,1][0,1,-1] \nonumber \\
 & \rightarrow & [0,0,0]
\end{eqnarray}
as expected.
\bigbreak
The reader will notice that it is unnecessary to keep track of the permutation induced by the first generator in a braid word (or the first braid word in a product, when the corresponding pure twist-word is known), as there is nothing for this permutation to act upon. 
\bigbreak 

\begin{figure}
     \begin{center}
     \begin{tabular}{c}
     \includegraphics[width=7cm]{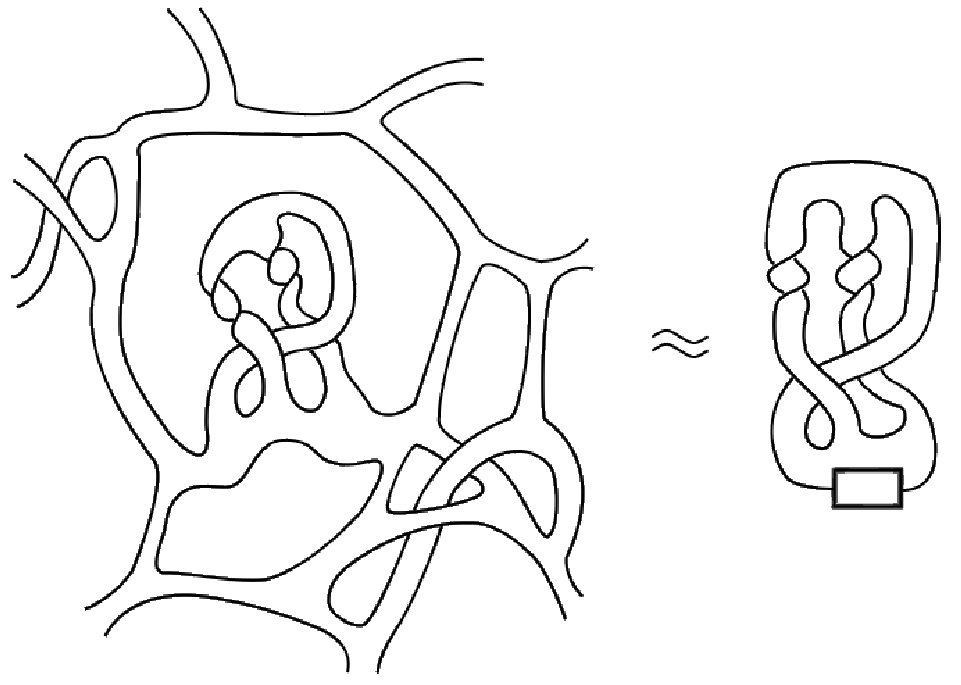}
     \end{tabular}
     \caption{A braided belt can be embedded within a larger network, represented by the box.}
     \label{fig:belt_in_network}
\end{center}
\end{figure}

The braided belts we consider in this paper are of interest not only from a purely topological basis, but also due to a possible connection with theoretical physics. In \cite{Sundance1} it was shown that braided belts attached at one end to a larger network of ribbons could be used to represent the elementary quarks and leptons. In order to keep the discussion in this paper relevant to the work in \cite{Sundance1}, we shall henceforth treat the ordering of generators in a braid word as indicative of a ``top end'' which is free to be trinion-flipped, and a ``bottom end'' which is attached to a larger network (which for all practical purposes is fixed and static). The left-most generator in a product is equated with the top end, and this is why we shall always work from left to right when we resolve a braid word to find the associated pure twist-word. In diagrams of braided belts we shall henceforth include a box on the boundary of the belt at the bottom end, to represent the presence of a larger network to which the braided belt is attached, as shown in Fig.~\ref{fig:belt_in_network} 
\bigbreak

\begin{figure}
     \begin{center}
     \begin{tabular}{c}
     \includegraphics[width=7cm]{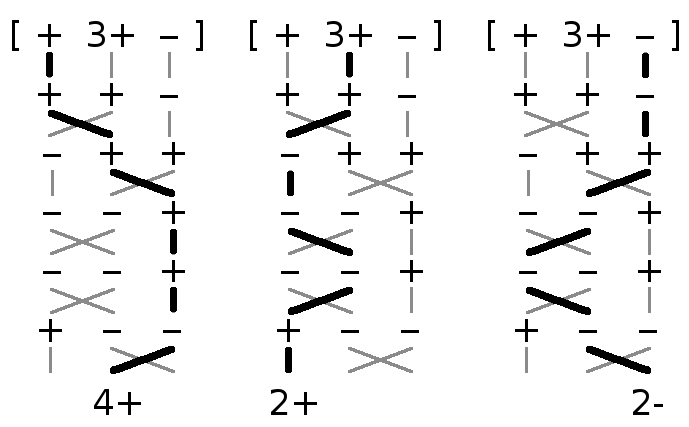}
     \end{tabular}
     \caption{Finding the pure twist-word for a braid $[\fh,\frac{3}{2},-\fh]\sigma_1\sigma_2\sigma^{-1}_1\sigma^{-1}_1\sigma^{-1}_2$}
     \label{fig:braidCalcul}
\end{center}
\end{figure}

There is a simple schematic technique for finding the pure twist-word associated with a braid word, as follows: 
\begin{itemize}
     \item Replace each generator in the braid word by one of the following triples of symbols;
	\begin{equation}
\begin{tabular}{lll}
$\sigma_1$      & $\rightarrow$ & $+\:+\:-$ \\
$\sigma_1^{-1}$ & $\rightarrow$ & $-\:-\:+$ \\
$\sigma_2$      & $\rightarrow$ & $-\:+\:+$ \\
$\sigma_2^{-1}$ & $\rightarrow$ & $+\:-\:-$
\end{tabular}
	\end{equation}
     such that they form a vertical stack, leftmost generator at the top.
     \item Beneath each pair of similar symbols, place a cross, $\times$, connecting each symbol to the diagonally opposite symbol in the pair below. Beneath each dissimilar symbol, place a vertical bar, $|$, connecting it directly to the symbol below it. The lines below the bottom row of symbols will end on three blank spaces.
     \item Sum up the symbols along each vertical path, starting at the top of the symbol stack, and ending on the blank spaces at the bottom. Each $+$ stands in for $\fh$. Each $-$ stands in for $-\fh$. The resulting triplet of numbers is the pure twist-word.
\end{itemize}
In the case where the braid is written in standard form $[a,b,c]{\cal B}$ we construct the symbol stack for the braid word as described above, and then place the twist word at the top, with three vertical lines descending to the top of the symbol stack. 

\bigbreak
In Fig.~\ref{fig:braidCalcul} we give an example of using this process to find the pure twist-word corresponding to the braid $[\fh,\frac{3}{2},-\fh]\sigma_1\sigma_2\sigma^{-1}_1\sigma^{-1}_1\sigma^{-1}_2$. The three diagrams in the figure correspond with the addition of twists along each of the three strands. The resulting pure twist-word is found to be $[1,2,-1]$.
\bigbreak

\begin{figure}
     \begin{center}
     \begin{tabular}{c}
     \includegraphics[width=7cm]{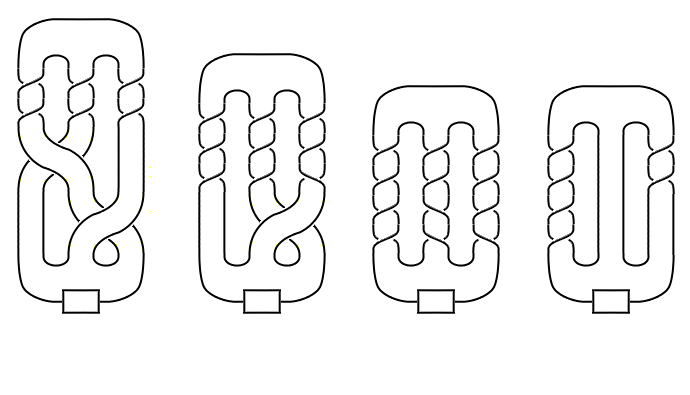}
     \end{tabular}
     \caption{The Positron, in fully-braided form (left), resolved in stages to its pure twist form (right).}
     \label{fig:positron}
     \end{center}
\end{figure}

A further example is given by the braid $e^{+}_\mathrm{L} = [1,1,1]\sigma_{1}\sigma_{2}^{-1}$. This is the braid structure assigned to the left-handed positron in \cite{Sundance}.
We find that 
\begin{eqnarray}
e^{+}_\mathrm{L} & \rightarrow & [1,1,1]\sigma_{1}\sigma_{2}^{-1} \nonumber \\
 & \rightarrow & [1,1,1]\left[\fh,\fh,-\fh\right]\sigma_{2}^{-1}  \nonumber \\
 & \rightarrow & \left[\frac{3}{2},\frac{3}{2},\fh\right]\sigma_{2}^{-1}  \nonumber \\
 & \rightarrow & P_{\sigma_{2}^{-1}}\left(\left[\frac{3}{2},\frac{3}{2},\fh\right]\right)
\left[\fh,-\fh,-\fh\right] \nonumber \\
 & \rightarrow & \left[\frac{3}{2},\fh,\frac{3}{2}\right]\left[\fh,-\fh,-\fh\right] \nonumber \\
 & \rightarrow & [2,0,1].
\end{eqnarray}
Thus $[2,0,1]$ gives the framings on the equivalent parallel flat strip belt. See Fig.~\ref{fig:positron} for a graphical version of this calculation, and a depiction of the 
boundary of the surface that corresponds to $e^{+}_\mathrm{L}.$ We denote the boundary of this surface by $\partial e^{+}_\mathrm{L}$. Note that $\partial e^{+}_\mathrm{L}$ is independent (topologically) of the deformation that we have applied to straighten out the braiding from
the original definition of $e^{+}_\mathrm{L}$. Thus our algebraic reduction gives us an algorithm for finding the boundary link for each particle in Bilson-Thompson's tables. \bigbreak

\begin{figure}
     \begin{center}
     \begin{tabular}{c}
     \includegraphics[width=7cm]{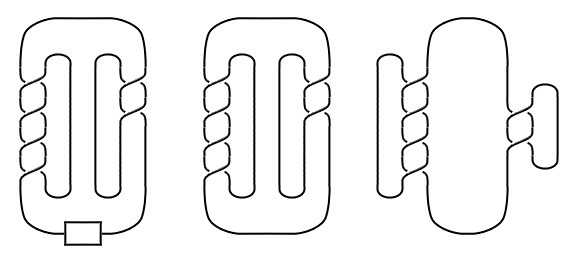}
     \end{tabular}
     \caption{By assuming that the outer edge of the braided belt is closed (when we trace it through the rest of the network), we can equate a link (right) to the pure twist form (left) of any braid, in this case the left-handed positron as illustrated in Fig~\ref{fig:positron}}
     \label{fig:positron_boundary}
\end{center}
\end{figure}

\begin{figure}
     \begin{center}
     \begin{tabular}{c}
     \includegraphics[width=7cm]{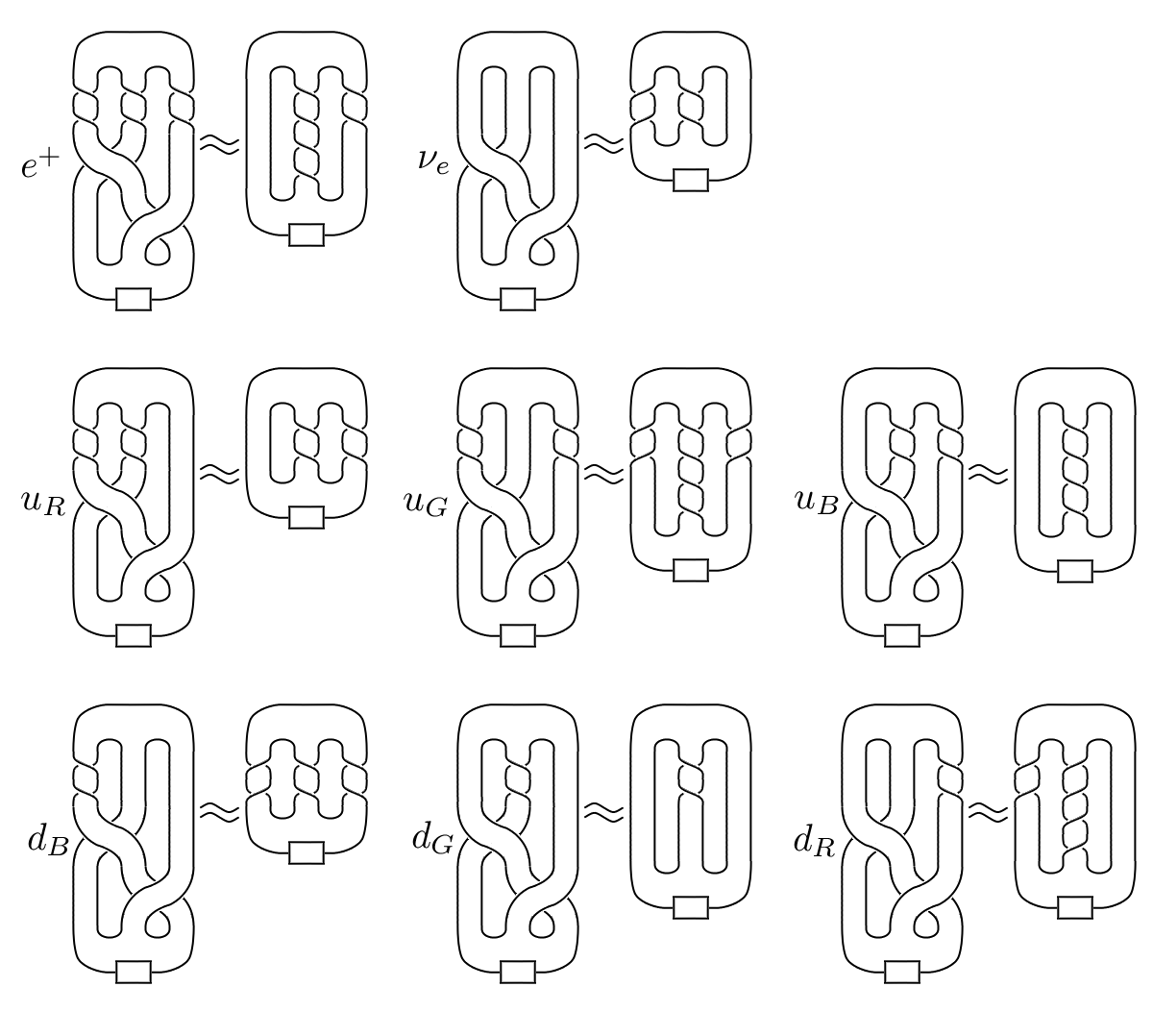}
     \end{tabular}
     \caption{Left-handed negatively-charged fermions, as per the structure proposed by Bilson-Thompson, and their associated pure twist form}
     \label{fig:L-neg}
\end{center}
\end{figure}

\begin{figure}
     \begin{center}
     \begin{tabular}{c}
     \includegraphics[width=7cm]{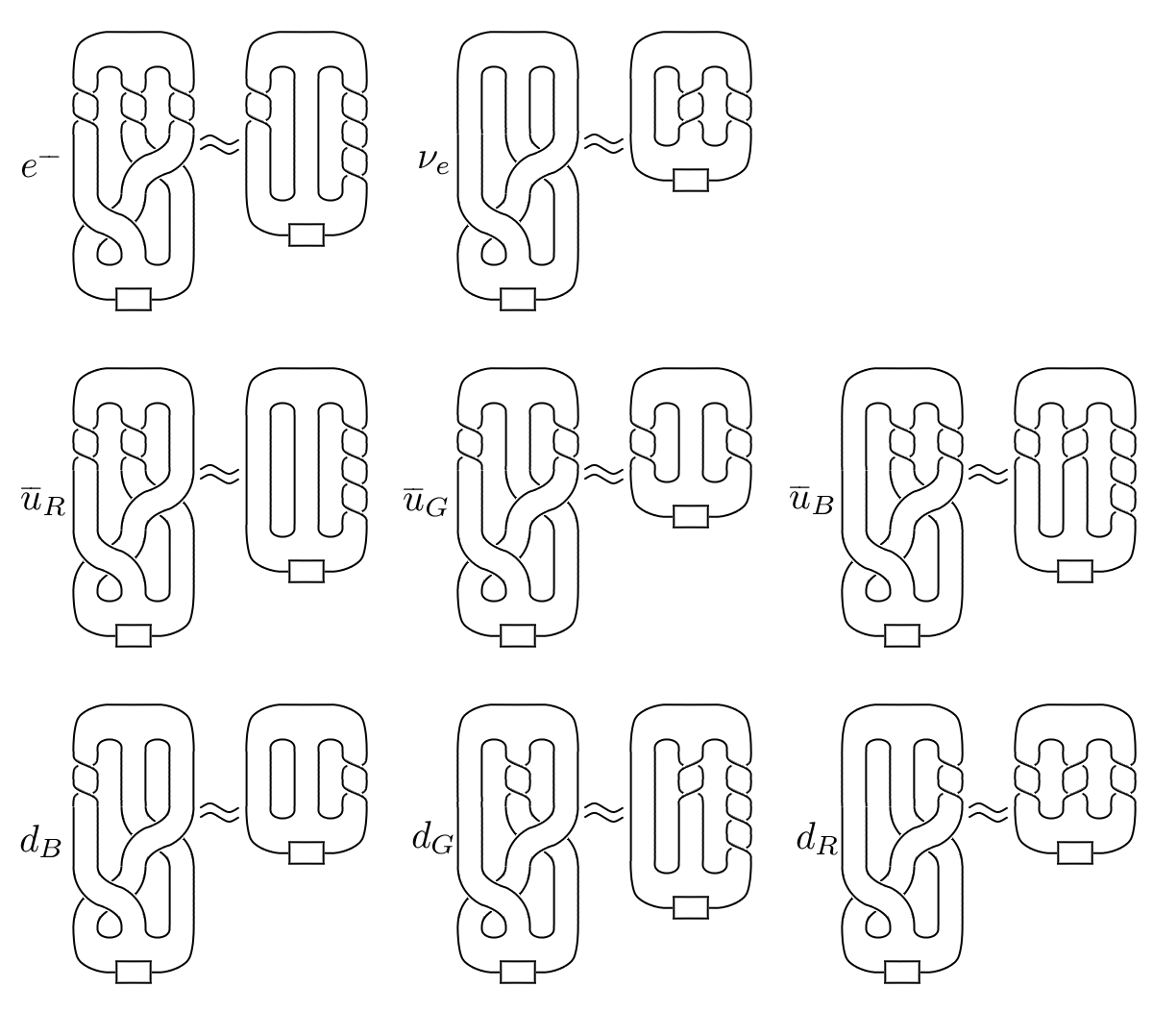}
     \end{tabular}
     \caption{Right-handed negatively-charged fermions, as per the structure proposed by Bilson-Thompson, and their associated pure twist form}
     \label{fig:R-neg}
\end{center}
\end{figure}

\begin{figure}
     \begin{center}
     \begin{tabular}{c}
     \includegraphics[width=7cm]{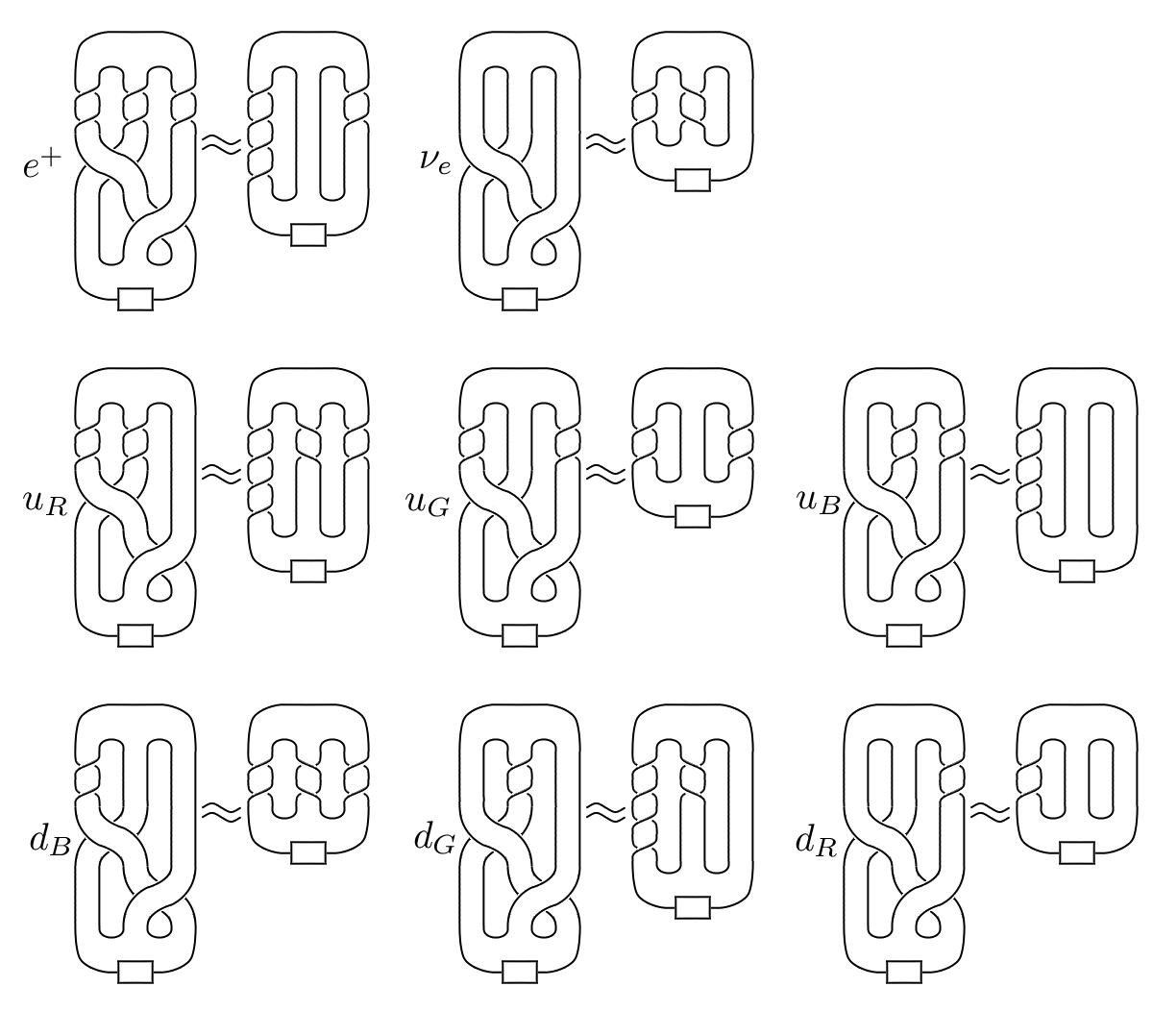}
     \end{tabular}
     \caption{Left-handed positively-charged fermions, as per the structure proposed by Bilson-Thompson, and their associated pure twist form}
     \label{fig:L-pos}
\end{center}
\end{figure}

\begin{figure}
     \begin{center}
     \begin{tabular}{c}
     \includegraphics[width=7cm]{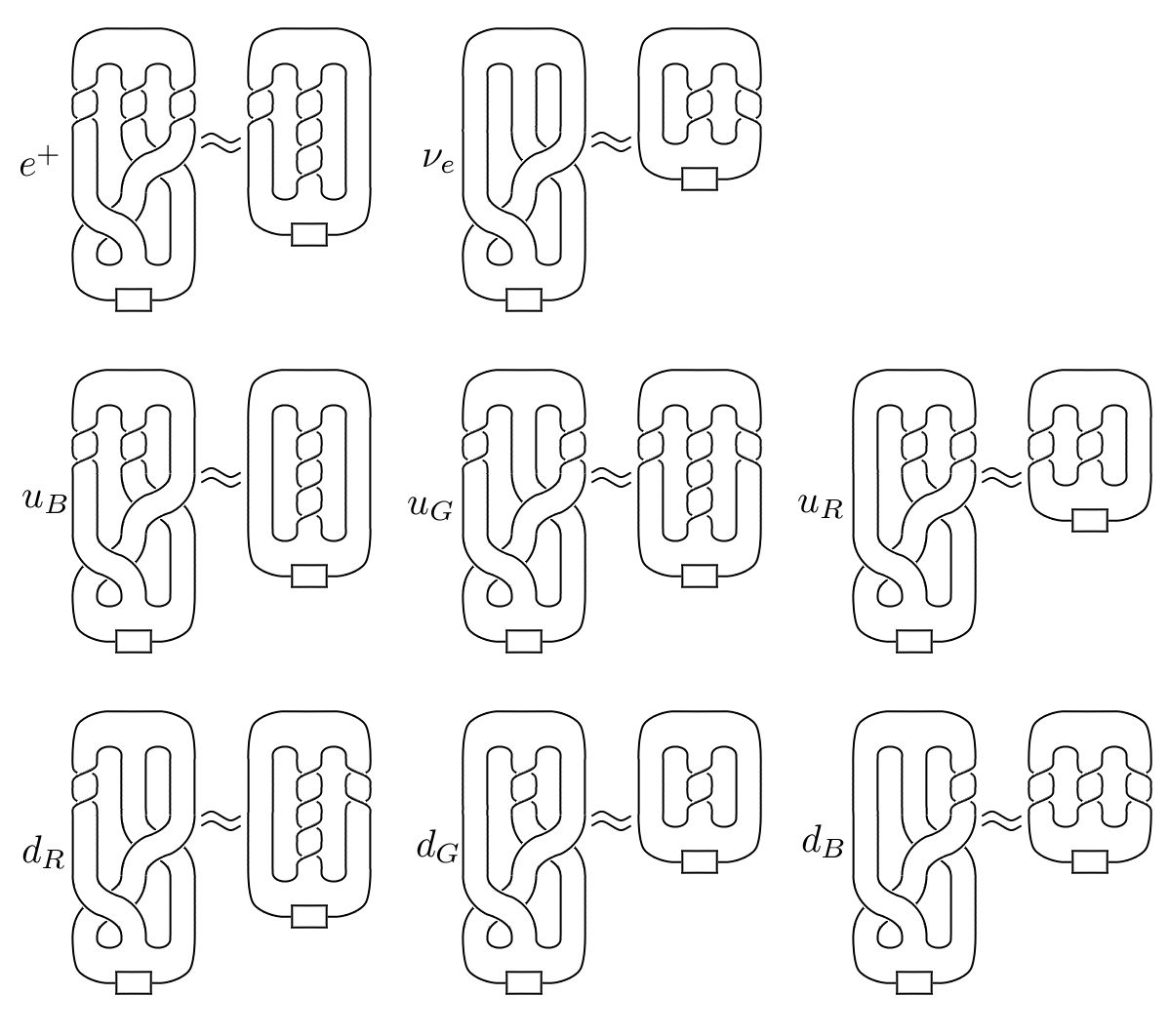}
     \end{tabular}
     \caption{Right-handed positively-charged fermions, as per the structure proposed by Bilson-Thompson, and their associated pure twist form}
     \label{fig:R-pos}
\end{center}
\end{figure}

In Fig.~\ref{fig:L-neg} - Fig.~\ref{fig:R-pos} we illustrate the correspondence between the braids proposed in \cite{Sundance} to match the first-generation fermions of the Standard Model, and their pure twist form. Notice that the pure twist forms are distinct in each illustrated case, except for the neutrino and anti-neutrino. In this case, the left-handed ``negative'' neutrino is isomorphic to the right-handed ``positive'' anti-neutrino. Likewise the left-handed ``positive'' anti-neutrino is isomorphic to the right-handed ``negative'' neutrino. In other words, there are half as many topologically distinct states for neutral particles as one would expect for charged particles. This is in agreement with the Standard Model, where neutrinos are purely left-handed, and anti-neutrinos are purely right-handed. The pure twist numbers in each case are listed with the corresponding particles in the table below;
\bigbreak

\begin{center}
\begin{tabular}{cc||cc}
left-handed &           &right-handed &   \\ \hline
$e^-$       & [0,-2,-1] & $e^-$       & [-1,0,-2] \\
$\bar{u}_B$ & [0,-1,-1] & $\bar{u}_B$ & [-1,1,-2] \\
$\bar{u}_G$ & [1,-2,-1] & $\bar{u}_G$ & [-1,0,-1] \\
$\bar{u}_R$ &  [0,-2,0] & $\bar{u}_R$ &  [0,0,-2] \\
$d_B$       &  [1,-2,0] & $d_B$       &  [0,0,-1] \\
$d_G$       &  [0,-1,0] & $d_G$       &  [0,1,-2] \\
$d_R$       & [1,-1,-1] & $d_R$       & [-1,1,-1] \\
$\nu_L$     &  [1,-1,0] & $\nu_R$     &  [0,1,-1] \\
$\bar{d}_B$ &   [1,0,0] & $\bar{d}_B$ &  [0,2,-1] \\
$\bar{d}_G$ &  [2,-1,0] & $\bar{d}_G$ &   [0,1,0] \\
$\bar{d}_R$ &  [1,-1,1] & $\bar{d}_R$ &  [1,1,-1] \\
$u_B$       &  [2,-1,1] & $u_B$       &   [1,1,0] \\
$u_G$       &   [1,0,1] & $u_G$       &  [1,2,-1] \\
$u_R$       &   [2,0,0] & $u_R$       &   [0,2,0] \\
$e^+$       &   [2,0,1] & $e^+$       &   [1,2,0]
\end{tabular}
\end{center}

\bigbreak

\section{Algebra}
\label{sec:algebra}
In this section we give a representation of the three-strand braid group in terms of {\it permutation matrices}. A permutation matrix $P$ is a matrix whose columns are a permutation 
of the columns of the identity matrix. Such a matrix acts as a permutation on the standard basis (column vectors that have a single unit entry). We shall sometimes refer to $P$ as a permutation, rather than a permutation matrix, for brevity. For example,
\beq
P_{1,2} =
\left( \begin{array}{ccc}
0 & 1 & 0\\
1 & 0 & 0 \\
0 & 0 & 1\\
\end{array} \right)
\eeq 
and
\beq
P_{2,3} =
\left( \begin{array}{ccc}
1 & 0 & 0\\
0 & 0 & 1 \\
0 & 1 & 0\\
\end{array} \right).
\eeq
Obviously 
$$\left( \begin{array}{ccc}
0 & 1 & 0\\
1 & 0 & 0 \\
0 & 0 & 1\\
\end{array} \right)\left(\begin{array}{c}1\\0\\0\end{array}\right)
= \left(\begin{array}{c}0\\1\\0\end{array}\right)$$
and so on.
\bigbreak

If $D$ is a diagonal matrix and $P$ is a permutation matrix, then 
\beq
PD = D^P P 
\label{eq:PDtoDP}
\eeq
where $D^P = P(D)$ denotes the result of permuting the elements of $D$ along the 
diagonal according to the permutation $P$. This is directly analogous to eqn.~(\ref{eq:wordpermute})
\bigbreak

The permutation matrices $P_{1,2}$ and $P_{2,3}$ generate (by taking products) all permutations on three letters. These permutations can be used to represent permutations induced by a braid, and hence they can stand in for the braid word. A general twist word 
$[a,b,c]$ may be represented by the matrix
$$[a,b,c] =
\left( \begin{array}{ccc}
t^{a} & 0 & 0\\
0 & t^{b} & 0 \\
0 & 0 & t^{c}\\
\end{array} \right).$$
The reader should be aware that the permutation matrices do not give a complete image of the braid group, because they do not define the direction of crossing. Hence the $P$s contain less information than the $\sigma_i$s. However we can find the pure twist form of a braid using this matrix representation in a manner we will now describe.
\bigbreak

Consider the identity 3-belt, as illustrated in fig.~\ref{fig:closed}. If we perform trinion flips on the top of the belt, we induce both crossings and twistings. When the crossings created in this process correspond to a given generator, $\sigma_i$, the twists induced are the {\em negative} of the twists corresponding to $\sigma_i$, as listed in eqns~(\ref{eq:sigmaToTwist}). The relative minus sign occurs because eqns~(\ref{eq:sigmaToTwist}) list the twist words created when crossings are eliminated, rather than created. In the interests of notational clarity, we shall therefore write $\rho_i$ to denote the twists {\em and} permutations induced on an initially trivial 3-belt by performing a trinion flip to cross strand $i$ over strand $i+1$. The corresponding twists and permutations are then; 
\begin{equation}
\begin{tabular}{lll}
$\rho_1$      & $\rightarrow$ & $\left[-\fh, \:-\fh,\:\fh\right]P_{1,2}$ \vspace{1mm} \\
$\rho_1^{-1}$ & $\rightarrow$ & $\left[\fh,\:\fh,\:-\fh\right]P_{1,2}$ \vspace{1mm} \\
$\rho_2$      & $\rightarrow$ & $\left[\fh,\:-\fh, \:-\fh\right]P_{2,3}$ \vspace{1mm} \\
$\rho_2^{-1}$ & $\rightarrow$ & $\left[-\fh,\:\fh,\:\fh\right]P_{2,3}$
\end{tabular}
\label{eq:rho}
\end{equation}
compare this with eqns~(\ref{eq:sigmaToTwist}). Applying eqn.~(\ref{eq:PDtoDP}) we see that
$$P_{1,2}[a,b,c] = [b,a,c]P_{1,2}$$
and
$$P_{2,3}[a,b,c] = [a,c,b]P_{2,3}.$$
\bigbreak

To find the pure twist-word corresponding to a given braid word, we first replace the $\sigma_i$ with the equivalent $\rho_i$ e.g.
$$
\sigma_1 \sigma_2 \sigma^{-1}_1 \rightarrow \rho_1 \rho_2 \rho^{-1}_1.
$$
We next substitute in the twist words and permutations from eqns~(\ref{eq:rho}), and apply eqn.~(\ref{eq:PDtoDP}). Once all the twist words have been shifted to the far left of the resulting expression and summed, we read off the negative of this twist word to find the pure twist-word corresponding to our initial braid word.
\bigbreak

\noindent{\bf Example:} Let us find the pure twist word corresponding to the braid $\sigma_1\sigma_2\sigma^{-1}_1$. We proceed as follows;
\begin{eqnarray}
\rho_1\rho_2\rho^{-1}_1 & = & \left[-\fh,\:-\fh,\:\fh\right]P_{1,2}\left[\fh,\:-\fh,\:-\fh\right] P_{2,3}\left[\fh,\:\fh,\:-\fh\right]P_{1,2} \nonumber \\
 & = & \left[-\fh,\:-\fh,\:\fh\right]\left[-\fh,\:\fh,\:-\fh\right]P_{1,2}P_{2,3} \left[\fh,\:\fh,\:-\fh\right]P_{1,2} \nonumber \\
 & = & \left[-\fh,\:-\fh,\:\fh\right]\left[-\fh,\:\fh,\:-\fh\right] \left[-\fh,\:\fh,\:\fh\right]P_{1,2}P_{2,3}P_{1,2} \nonumber \\
 & = & \left[-\frac{3}{2},\:\fh,\:\fh\right]P_{1,2}P_{2,3}P_{1,2} \nonumber 
\end{eqnarray}
We then take the negative of the computed twist word, to obtain our result, $[\frac{3}{2},-\fh,-\fh]$.
\bigbreak

These assignments of framed permutations to braids gives a representation of the framed braid group into framed permutations. To see this we can directly verify that $\sigma_{1}\sigma_{2}\sigma_{1} = \sigma_{2}\sigma_{1}\sigma_{2}$
at the level of the framed permutations as follows; 
\begin{eqnarray}
\rho_1 \rho_2 \rho_1 & = & \left[-\fh,-\fh,\fh\right]P_{1,2}\left[\fh,-\fh,-\fh\right]P_{2,3}\left[-\fh,-\fh,\fh\right]P_{1,2} \nonumber \\
 & = & \left[-\fh,-\fh,\fh\right]\left[-\fh,\fh,-\fh\right]P_{1,2}P_{2,3}\left[-\fh,-\fh,\fh\right]P_{1,2} \nonumber \\
 & = & \left[-1,0,0\right]P_{1,2}P_{2,3}\left[-\fh,-\fh,\fh\right]P_{1,2} \nonumber \\
 & = & \left[-1,0,0\right]P_{1,2}\left[-\fh,\fh,-\fh\right]P_{2,3}P_{1,2} \nonumber \\
 & = & \left[-1,0,0\right]\left[\fh,-\fh,-\fh\right]P_{1,2}P_{2,3}P_{1,2} \nonumber \\
 & = & \left[-\fh,-\fh,-\fh\right]P_{1,2}P_{2,3}P_{1,2}.
\end{eqnarray}
Similarly, we find that
$$\rho_2 \rho_1 \rho_2 = \left[-\fh,-\fh,-\fh\right]P_{2,3}P_{1,2}P_{2,3}.$$ 
Since $\rho_{1} \rho_{2} \rho_{1}$ and $\rho_{2} \rho_{1}\rho_{2}$ yield the same twist word we conclude that $\sigma_{1} \sigma_{2} \sigma_{1} = \sigma_{2} \sigma_{1}\sigma_{2}.$  
It is noteworthy that $P_{1,2}P_{2,3}P_{1,2} = P_{2,3}P_{1,2}P_{2,3}$, however this has no bearing on the result because the twist words alone are sufficient to define isomorphism.
\bigbreak

\begin{figure}
     \begin{center}
     \begin{tabular}{c}
     \includegraphics[width=6cm]{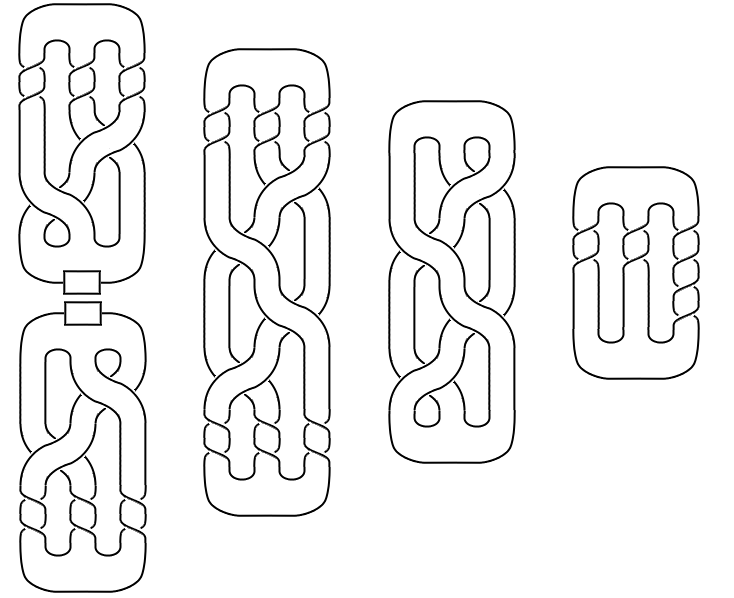}
     \end{tabular}
     \caption{Uniquely forming the product of two braids (left) by joining them box-to-box, to yield a link (right).}
     \label{fig:combining}
\end{center}
\end{figure}

\section{Conclusions}
\label{sec:concl}
We have demonstrated several equivalent approaches to manipulating framed braids on three strands, carried on a surface with boundary (i.e. ``braided belts''), which yield an invariant form having no crossings. The twists on each strand in this form define a triplet of numbers. Any braided belts having the same triplet of associated numbers (the pure twist word) are isotopically equivalent. We have furthermore demonstrated that the 3-braids proposed in \cite{Sundance} to represent the fermions of the first generation are indeed isotopically inequivalent.
\bigbreak

\begin{figure}
     \begin{center}
     \begin{tabular}{c}
     \includegraphics[width=9cm]{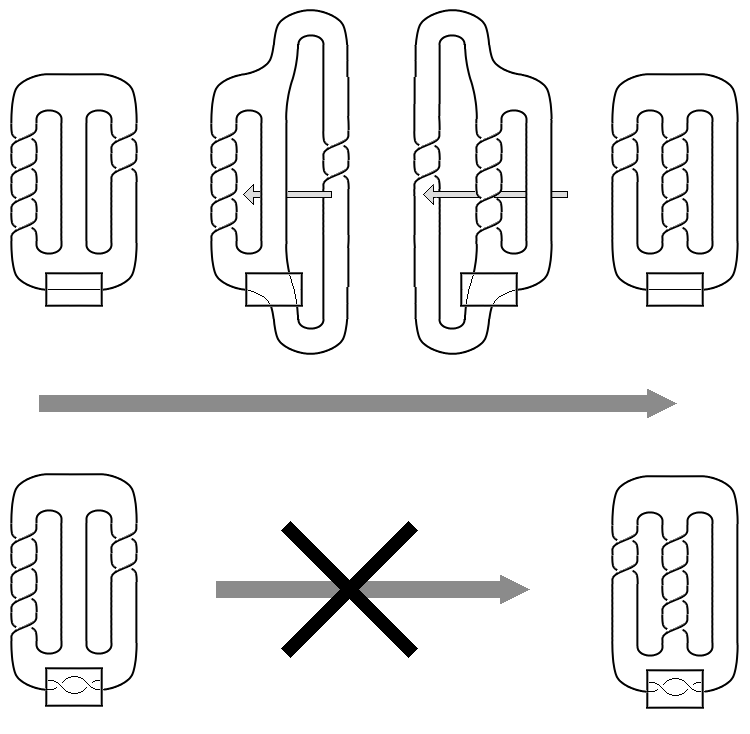}
     \end{tabular}
     \caption{If the edge running through the external network (represented by the box) carries no twists, other than those on the braided belt strands themselves, the strands may be cyclically permuted by passing one strand under (or around) the others (top row, left to right). Fig.~\ref{fig:belt_in_network} gives an explicit example of such a case. If the edge does twist around other edges in the rest of the network (bottom row) then such cyclic permutations are not allowed.}
     \label{fig:reordering}
\end{center}
\end{figure}

We have demonstrated how a braided belt on three strands may be mapped to a link, with the linking numbers being computable from the pure twist numbers. We should note that the usual ambiguity then arises, as to how we uniquely define a product of these links (which would be equivalent to interacting two particles with each other), since the links may be deformed by Reidemeister moves, leading to different possible products (depending on the point at which the links are joined to each other). In \cite{Sundance}, the braids were assumed to have a top and bottom end, and particles would annihilate with their corresponding antiparticles by joining ``bottom-to-bottom''. We may adapt this approach by treating the box (representing the rest of the network within which the braids are embedded) as the bottom of any braid, and combining braids (or links) box-to-box to form their product, as illustrated in Fig.~\ref{fig:combining}. For the moment this may be treated as a convenience. The physical interpretation of this requirement will be left for future work. 
\bigbreak

When considering a braided belt embedded at one end in a network, as we have discussed in this paper, the pure twist-word of the braid is only defined up to cyclic permutations of $[a,b,c]$. A single permutation is uniquely specified only if the path along the outside edge of the braid and through the network has topology that obstructs this cyclic permutation. In many cases this will happen when the network contains nontrivial twists or linking. This may be seen by recognising that the strands exiting a node can in general be cyclically permuted around the node. If the external edge of the network is sufficiently simple, the network
itself may be treated as a node (Fig.~\ref{fig:reordering}). This clearly has consequences for the interpretation of particle states proposed in \cite{Sundance}, and will also be addressed in future work.\bigbreak

The focus of this paper has been on manipulation of braids, rather than direct connections to physics. However this work has established a useful basis for discussions of physics in future work. In particular, functions of the pure twist words of various fermions can be matched to quantum numbers like weak isospin and hypercharge, baryon number, lepton number, and so forth.


\begin{thebibliography}{99} 

\bibitem{Sundance}
Sundance O. Bilson-Thompson, A topological model of composite preons, arXiv:hep-th/0503213

\bibitem{Sundance1}
Sundance O. Bilson-Thompson, Fotini Markopoulou, Lee Smolin, Quantum gravity and the standard model, 
{\it Class. Quant. Gravity} {\bf 24} (2007) 3975-3993, arXiv:hep-th/0603022

\bibitem{Hackett2}
Sundance Bilson-Thompson, Jonathan Hackett, Louis H. Kauffman, Lee Smolin,
Particle Identifications from Symmetries of Braided Ribbon Network Invariants,
aXiv:0804.0037
 
\bibitem{Sossinsky}
V. V. Prasolov and A. B. Sossinsky, {\em Knots, Links, Braids and 3-Manifolds,} American Mathematical Society - 
Translations of Mathematical Monographs, Volume 154 (1997).

\bibitem{AT} 
M.F.\ Atiyah, {\em The Geometry and Physics of Knots},  Cambridge University
Press, 1990.

\bibitem{Birman}
J. Birman, {\em Braids, Links and Mapping Class Groups,} {\em Ann. of Math. Studies No. 82}, Princeton, N.J.: Princeton University
Press (1976).

\bibitem{Conway}
J. H. Conway, An enumeration of knots and links and some of their related properties, In: ``Computational Problems in Abstract Algebra" (Oxford 1970), 329 - 358, Pergamon Press.


\bibitem {KP}
L.H. Kauffman, {\em Knots and Physics}, World Scientific Publishers (1991), 
Second Edition (1993), Third Edition (2002).


\bibitem{WIT} 
E. Witten, Quantum field theory and the Jones polynomial,
Commun. Math. Phys.  {\bf 121} (1989), 351--399.


 
\end{thebibliography}
\end{document}